\documentclass[12pt]{article}
\usepackage{amssymb,amsthm,amsmath,amsfonts,latexsym,tikz,hyperref}
\usepackage[hmargin=.8in,vmargin=.5in]{geometry}

\newtheorem{thm}{Theorem}[section]
\newtheorem{prop}[thm]{Proposition}
\newtheorem{cor}[thm]{Corollary}
\newtheorem{lem}[thm]{Lemma}
\newtheorem{conj}[thm]{Conjecture}
\newtheorem{exa}[thm]{Example}

\newcommand{\lf}{\lfloor}
\newcommand{\rf}{\rfloor}

\newcommand{\lc}{\lceil}
\newcommand{\rc}{\rceil}

\newcommand{\bin}[2]{\left\{ #1 \atop #2 \right\}}

\newcommand{\ben}{\begin{enumerate}}
\newcommand{\een}{\end{enumerate}}
\newcommand{\ble}{\begin{lem}}
\newcommand{\ele}{\end{lem}}
\newcommand{\bth}{\begin{thm}}
\renewcommand{\eth}{\end{thm}}
\newcommand{\bpr}{\begin{prop}}
\newcommand{\epr}{\end{prop}}
\newcommand{\bco}{\begin{cor}}
\newcommand{\eco}{\end{cor}}
\newcommand{\bcon}{\begin{conj}}
\newcommand{\econ}{\end{conj}}
\newcommand{\bde}{\begin{defn}}
\newcommand{\ede}{\end{defn}}
\newcommand{\bex}{\begin{exa}}
\newcommand{\eex}{\end{exa}}
\newcommand{\barr}{\begin{array}}
\newcommand{\earr}{\end{array}}
\newcommand{\btab}{\begin{tabular}}
\newcommand{\etab}{\end{tabular}}
\newcommand{\beq}{\begin{equation}}
\newcommand{\eeq}{\end{equation}}
\newcommand{\bea}{\begin{eqnarray*}}
\newcommand{\eea}{\end{eqnarray*}}
\newcommand{\bal}{\begin{align*}}
\newcommand{\bce}{\begin{center}}
\newcommand{\ece}{\end{center}}
\newcommand{\bpi}{\begin{picture}}
\newcommand{\epi}{\end{picture}}
\newcommand{\bpp}{\begin{picture}}
\newcommand{\epp}{\end{picture}}
\newcommand{\bfi}{\begin{figure} \begin{center}}
\newcommand{\efi}{\end{center} \end{figure}}
\newcommand{\bprf}{\begin{proof}}
\newcommand{\eprf}{\end{proof}\medskip}
\newcommand{\capt}{\caption}
\newcommand{\bsl}{\begin{slide}{}}
\newcommand{\esl}{\end{slide}}
\newcommand{\bfr}{\begin{frame}}
\newcommand{\efr}{\end{frame}}

\newcommand{\hqed}{\hfill \qed}
\newcommand{\hqedm}{\hfill \qed \medskip}

\newcommand{\eqqed}[1]{$\rule{1ex}{0ex}\hfill{\dil#1}\hfill\qed$}

\newcommand{\hs}[1]{\hspace{#1}}
\newcommand{\hso}[1]{\hspace{-1pt}}
\newcommand{\vs}[1]{\vspace{#1}}

\newcommand{\emp}{\emptyset}

\newcommand{\sbe}{\subseteq}

\newcommand{\Cong}{\equiv}

\newcommand{\case}[4]{\left\{\barr{ll}#1&\mbox{#2}\\#3&\mbox{#4}\earr\right.}

\def\<{\langle}
\def\>{\rangle}

\newcommand{\ree}[1]{(\ref{#1})}

\newcommand{\llra}{\longleftrightarrow}

\newcommand{\io}{\iota}

\newcommand{\la}{\lambda}

\newcommand{\ze}{\zeta}

\newcommand{\bb}{{\bf b}}

\newcommand{\bF}{{\bf F}}

\newcommand{\bT}{{\bf T}}

\newcommand{\bbN}{{\mathbb N}}

\newcommand{\cD}{{\cal D}}

\newcommand{\cT}{{\cal T}}

\DeclareMathOperator{\Mod}{mod}

\newcommand{\dil}{\displaystyle}

\begin{document}
\pagestyle{plain}

\title{The fractal nature of the  Fibonomial triangle
}
\author{
Xi Chen\\[-5pt]
\small School of Mathematical Sciences, Dalian University of Technology,\\[-5pt]
\small Dalian City, Liaoning Province, 116024, P. R. China, {\tt xichen.dut@gmail.com} 
\thanks{Research partially supported by a grant from the China Scholarship Council}\\
and\\
Bruce E. Sagan\\[-5pt]
\small Department of Mathematics, Michigan State University,\\[-5pt]
\small East Lansing, MI 48824-1027, USA, {\tt sagan@math.msu.edu}
}

\date{\today\\[10pt]
	\begin{flushleft}
	\small Key Words: binomial coefficient, fibonomial coefficient, fractal, Lucas' congruence, tiling
	                                       \\[5pt]
	\small AMS subject classification (2000): 
	Primary 05A10;
	Secondary 11B39, 11B50, 11B65.
	\end{flushleft}}

\maketitle

\begin{abstract}
It is well known that Pascal's triangle exhibits fractal behavior when reduced modulo a prime.  We show that the triangle of Fibonomial coefficients has a similar nature modulo two.  Specifically, for any $m\ge0$, the subtriangle consisting of the first $3\cdot2^m$ rows is duplicated on the left and right sides of the next $3\cdot 2^m$ rows, with an inverted triangle of zeros in between.  We give three proofs of this fact.  The first uses a combinatorial interpretation of the Fibonomials due to Sagan and Savage.  The second employs an analogue of Lucas' congruence for the parity of binomial coefficients.  The final one is inductive.  We also use induction to show that the Fibonomial triangle has a similar structure modulo three.  We end with some open questions.
\end{abstract}

%
%

\section{Introduction}

Let $p$ be a prime and consider Pascal's triangle modulo $p$, $T_p$.  This array has fractal properties; see, for example, the article of Sved~\cite{sve:dv}.  We will describe this behavior more precisely for $p=2$ as the result will be needed in the sequel.  We will use $\Cong_2$ to denote congruence modulo two.   Let $\bbN$ denote the nonnegative integers and $m\in\bbN$.  Consider the triangle $T$ consisting of the first $2^m$ rows of $T_2$.  Then the next $2^m$ rows consist of a copy of $T$ against the left border, a copy of $T$ against the right border, and an inverted triangle of zeros between the two copies, see Figure~\ref{triangles}.  In terms of congruences, if 
$0\le n<2^m$ and $0\le k \le n+2^m$ then
$$
\binom{n+2^m}{k}\Cong_2
\begin{cases}
\dil\binom{n}{k}&\text{if $0\le k\le n$,}\\[10pt]
0 &\text{if $n<k<2^m$,}\\[10pt]
\dil\binom{n}{k-2^m}&\text{if $2^m\le k\le n+2^m$.}
\end{cases}
$$
These three cases can be summarized by the single formula
$$
\binom{n+2^m}{k}\Cong_2 \binom{n}{k}\qquad \text{if $0\le k<2^m$.}
$$
To see why the second case follows just note that $\binom{n}{k}=0$ for $k<0$ and $k>n$.  And the third case is obtained  from the previous equation and the symmetry of the binomial coefficients because, when $2^m\le k\le n+2^m$,
$$
\binom{n+2^m}{k}=\binom{n+2^m}{n+2^m-k}\Cong_2 \binom{n}{n+2^m-k}=\binom{n}{k-2^m}.
$$
We record this result for future use.
\bth
\label{binom}
Given $m\ge0$ and $0\le n,k<2^m$ we have

\medskip

\eqqed{
\binom{n+2^m}{k}\Cong_2 \binom{n}{k}.
}
\eth

\bfi
$$
\barr{ccccccccccccccc}
&&&&&&&1&&&&&&&\\
&&&&&&1&&1&&&&&&\\
&&&&&1&&0&&1&&&&&\\
&&&&1&&1&&1&&1&&&&\\
&&&1&&0&&0&&0&&1&&&\\
&&1&&1&&0&&0&&1&&1&&\\
&1&&0&&1&&0&&1&&0&&1&\\
1&&1&&1&&1&&1&&1&&1&&1
\earr
\hs{40pt}
\barr{ccccccccccc}
&&&&&1&&&&&\\
&&&&1&&1&&&&\\
&&&1&&1&&1&&&\\
&&1&&0&&0&&1&&\\
&1&&1&&0&&1&&1&\\
1&&1&&1&&1&&1&&1\\
&&&&&&&&&&\\
&&&&&&&&&&\\
\earr
$$
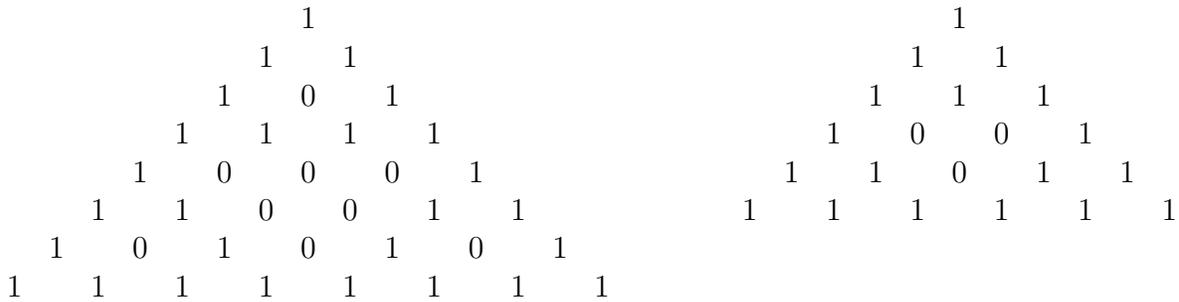
\capt{Pascal's triangle (left) and the Fibonomial triangle (right) modulo $2$}
\label{triangles} 
\efi

We now turn to the Fibonomial coefficients.  Consider the Fibonacci numbers as defined by $F_0=0$, $F_1=1$ and $F_n=F_{n-1}+F_{n-2}$ for $n\ge2$.  One then defines \emph{Fibotorials} and \emph{Fibonomial coefficients} by
$$
n!_F = F_1 F_2\cdots F_n
$$
and
$$
\binom{n}{k}_F=\frac{n!_F}{k!_F (n-k)!_F},
$$
respectively, where $0\le k\le n$.  It is convenient to let $\binom{n}{k}_F=0$ for $k<0$ and $k>n$.  Also note that the Fibonomials obviously display the same symmetry as the binomials in that $\binom{n}{k}_F=\binom{n}{n-k}_F$.  The purpose of this paper is to prove that the Fibonomial triangle displays a similar fractal behavior to Pascal's triangle but with triangles which have side length equal to $3\cdot 2^m$ for some $m\ge0$. Again, see Figure~\ref{triangles}.  Precisely stated,  our main theorem (for which we will give three proofs)  is as follows.
\bth
\label{main}
Given $m\ge0$ and $0\le n,k<3\cdot 2^m$ we have
$$
\binom{n+3\cdot 2^m}{k}_F\Cong_2 \binom{n}{k}_F.
$$
\eth

The rest of this paper is organized as follows.  In the next section we give an enumerative proof of Theorem~\ref{main} using a combinatorial interpretation of the Fibonomial coefficients due to Sagan and Savage~\cite{ss:cib}.    Interestingly, we also need Theorem~\ref{binom} for this demonstration.  Section~\ref{acl} contains a proof of our main theorem based on a Fibonomial analogue of Lucas' congruence for binomial coefficients~\cite{luc:cne}.  To prove this congruence, we use a result of Knuth and Wilf~\cite{kw:ppd} giving the $p$-adic valuation of generalized binomial coefficients.  In Section~\ref{ip} we give our third, inductive proof of Theorem~\ref{main}.    Section~\ref{mtc} is devoted to a discussion of the Fibonomial triangle modulo three.  We indicate how our inductive proof for $p=2$ can be modified to cover this case.  We end with a section of open questions raised by this work.

%
%

\section{A combinatorial proof}
\label{cp}

We now give a combinatorial proof of Theorem~\ref{main}.  We start by reviewing one of the standard combinatorial interpretations of the Fibonacci numbers themselves.

A \emph{tiling}, $T$, of a row of squares is a covering of the squares with disjoint dominos (covering two squares) and monominos (covering one square).  We let 
$$
\cT_n =\{T\ :\ \text{$T$ is a tiling of a row of $n$ squares}\}.
$$
The elements of $\cT_3$ are given in Figure~\ref{T3}.   We consider the squares as labeled $1,\ldots,n$ from left to right. 

\bfi
\begin{tikzpicture}
\draw (0,0) grid (3,1);
\fill (.5,.5) circle (.1);
\fill (1.5,.5) circle (.1);
\fill (2.5,.5) circle (.1);
\end{tikzpicture}
\qquad
\begin{tikzpicture}
\draw (0,0) grid (3,1);
\draw (.5,.5)--(1.5,.5);
\fill (.5,.5) circle (.1);
\fill (1.5,.5) circle (.1);
\fill (2.5,.5) circle (.1);
\end{tikzpicture}
\qquad
\begin{tikzpicture}
\draw (0,0) grid (3,1);
\draw (1.5,.5)--(2.5,.5);
\fill (.5,.5) circle (.1);
\fill (1.5,.5) circle (.1);
\fill (2.5,.5) circle (.1);
\end{tikzpicture}
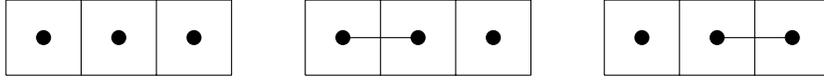
\capt{  The tilings in $\cT_3$}
\label{T3}
\efi

The next result is well known and follows directly from the definition of $\cT_n$ and the recursion for $F_n$.  Keep in mind the difference between the two subscripts for the future.
\ble
\label{Tn}
For $n\ge 1$ we have
$$
F_n = |\cT_{n-1}|
$$
where $|\cdot|$ denotes cardinality.\hqed
\ele

We must first investigate the behaviour of $F_n$ modulo two.  While the following proposition also follows easily by using the recursion, we will give a combinatorial proof in keeping with the focus of this section and because it will be useful in the sequel.
\bpr
\label{Fn2}
We have
$$
F_n \Cong_2 \case{0}{if $n \Cong_3 0$,}{1}{if $n\Cong_3 1$ or $2$.}
$$
\epr
\bprf
By Lemma~\ref{Tn}, it suffices to construct an involution $\io$ on $\cT_n$ which has no fixed points if $n\Cong_3 2$ and has exactly one fixed point otherwise.  We define $\io$ inductively as follows.

Consider the first two squares of $T\in\cT_n$.  If these squares are filled with two monominos, then replace them with a domino and vice-versa, keeping all other tiles the same.  This pairs up all tilings beginning with two monominos or with a domino.  The remaining unpaired $T$ must all begin with a monomino followed by a domino.  Now consider squares 4 and 5 of such a $T$.  Again, the tilings where those two squares contain two monominos or a domino are paired.  This is illustrated by the first diagram in Figure~\ref{io}.  This process is continued until only $r$ squares remain, where $r$ is the remainder of $n$ on division by $3$.  If $r=0$ or $1$ then we are left with a single fixed tiling.  The second drawing in Figure~\ref{io} illustrates this.  If $r=2$ then the last two squares must contain two monominos or a domino and so the last two tilings are paired leaving no fixed points.
\eprf

\bfi
\begin{tikzpicture}
\draw (0,0) grid (8,1);
\foreach \x  in {.5,1.5,...,7.5} \fill (\x,.5) circle (.1);
\draw (1.5,.5)--(2.5,.5) (6.5,.5)--(7.5,.5);
\end{tikzpicture}
\raisebox{12pt}{$\barr{c} \io \\[-5pt] \llra \earr$}
\begin{tikzpicture}
\draw (0,0) grid (8,1);
\foreach \x  in {.5,1.5,...,7.5} \fill (\x,.5) circle (.1);
\draw (1.5,.5)--(2.5,.5)  (3.5,.5)--(4.5,.5) (6.5,.5)--(7.5,.5);
\end{tikzpicture}

\vs{20pt}

\begin{tikzpicture}
\draw (0,0) grid (7,1);
\foreach \x  in {.5,1.5,...,6.5} \fill (\x,.5) circle (.1);
\draw (1.5,.5)--(2.5,.5) (4.5,.5)--(5.5,.5);
\pgfpathmoveto{\pgfpoint{7.5cm}{1cm}}
\pgfpatharc{90}{-90}{1cm and .5cm}
\pgfusepath{draw}
\draw[->] (7.5,1) -- (7.49,1);
\draw (8.8,.5) node {$\io$};
\end{tikzpicture}
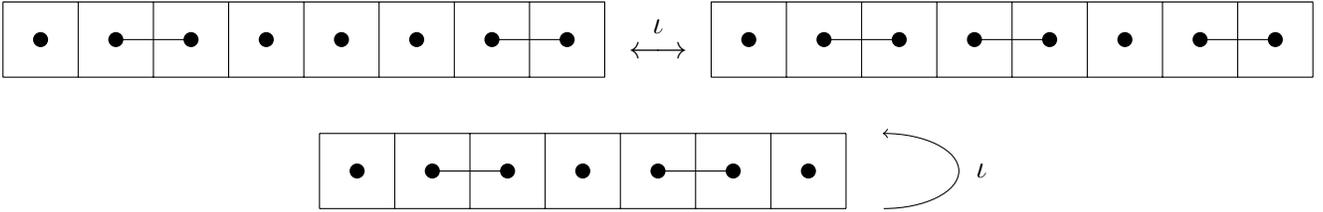
\capt{ A  pair (above) and fixed point (below) of $\io$}
\label{io} 
\efi

We are now in a position to recall the combinatorial interpretation of $\binom{n}{k}_F$ discovered by Sagan and Savage~\cite{ss:cib}.    An \emph{(integer) partition} is a weakly decreasing sequence $\la=(\la_1,\ldots,\la_l)$ of positive integers called \emph{parts}.  The associated \emph{Ferrers diagram} is an array of squares with $l$ left-justified rows and $\la_i$ squares in row $i$.  We write our Ferrers diagrams in English notation with the largest row on top and do not distinguish $\la$ from its Ferrers diagram.  The squares to the northwest of the heavy line in Figure~\ref{Tla} are the Ferrers diagram for the partition $(3,2,2,2)$.

A \emph{tiling} of $\la$ is a tiling of the parts of its Ferrers diagram, i.e., an element of $\cT_{\la_1}\times\cdots\times\cT_{\la_l}$.  We let
$$
\cT_\la=\{T\ :\ \text{$T$ is a tiling of $\la$}\}.
$$
The part of Figure~\ref{Tla} corresponding to $\la$ shows such a tiling.  We will also be interested in a subset of $\cT_\la$:
$$
\cD_\la=\{T\in \cT_\la\ |\ \text{each $T_{\la_i}$ starts with a domino}\}.
$$
Note that if $\la$ has a part of size one then $\cD_\la=\emp$.

We say that $\la$ \emph{fits in a $c\times d$ rectangle}, $\la\sbe c\times d$, if  the number of rows $\la$ is at most $c$ and the number of columns at most $d$.  So $\la$ fits in the upper left corner of the Ferrers diagram $c\times d$ consisting of $c$ rows of length $d$.  The squares in Figure~\ref{Tla} show the partition $(3,2,2,2)$ as it fits in a $4\times 5$ rectangle.  In this case, there is a \emph{complementary partition}, $\la^*$, whose parts are the columns of  the set difference $(c\times d)-\la$.  In Figure~\ref{Tla} we have $\la^*=(4,4,3)$.  Since the parts of $\la^*$ are columns, dominos in tilings of $\la^*$ will be vertical.  Figure~\ref{Tla} shows an element of $\cT_\la\times\cD_{\la^*}$ which  are the objects we need.
\bth[Sagan and Savage~\cite{ss:cib}]
\label{ss}
We have

\medskip

\eqqed{
\binom{n}{k}_F = \left\lvert \bigcup_{\la\sbe k\times (n-k)}   \cT_\la\times \cD_{\la^*}\right\rvert.
}

\eth

We can extend the action of the involution $\io$ in the proof of Proposition~\ref{Fn2} to $\cT_\la\times\cD_{\la^*}$ as follows.
Suppose $\la=(\la_1,\ldots,\la_l)$.  Apply $\io$ to the first row of $\la$, leaving the rest of the rectangle fixed.  If $\la_1\Cong_3 2$ then this will pair up all elements of $\cT_\la\times\cD_{\la^*}$.  If not, then the unpaired tilings are exactly those whose first row is the unique fixed point of $\io$ acting on $\cT_{\la_1}$.  Now iterate the process by applying $\io$ to the second row and subsequent rows if necessary.  If after finishing the rows of $\la$ we still have unpaired tilings, then we start with the columns of $\la^*$.  We first apply $\io$ to the portion of $\la_1^*$ above the leading domino and then continue until we have either paired up all tilings or have exactly one fixed point.  From  Proposition~\ref{Fn2} it is immediate that we will have a fixed point if and only if
\beq
\label{fp}
\text{$\la_i\Cong_3 0$ or $1$ for all $i$ and $\la_j^*\Cong_3 0$ or $2$ for all $j$.}
\eeq
This is a key step in the proof of the next result.

\bfi
\begin{tikzpicture}
\draw (0,0) grid (5,4);
\draw[very  thick] (0,0) -- (2,0) -- (2,3) -- (3,3) -- (3,4) -- (5,4);
\draw (.5,1.5) -- (1.5,1.5) (1.5,3.5) -- (2.5,3.5) (2.5,.5) -- (2.5,1.5) (3.5,.5) -- (3.5,1.5) (3.5,2.5) -- (3.5,3.5) (4.5,.5) -- (4.5,1.5);
\foreach \x  in {.5,1.5,...,4.5} \foreach \y in {.5,1.5,2.5,3.5}\fill (\x,\y) circle (.1);
\end{tikzpicture}
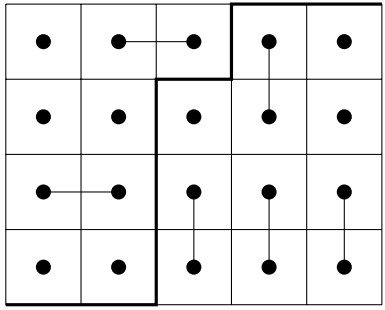
\capt{A tiling in $\cT_\la\times\cD_{\la^*}$ for $\la=(3,2,2,2)\sbe 4\times 5$}
\label{Tla} 
\efi

\ble
\label{lattice}
We have
$$
\binom{n}{k}_F\Cong_2
\begin{cases}
0 & \text{if $n\Cong_3 0$ and $k\Cong_3 1$,}\\[5pt]
\dil \binom{\lc 2n/3 \rc}{\lc 2k/3 \rc} & \text{if $n\Cong_3 1$ and $k\Cong_3 0$,}\\[20pt]
\dil  \binom{\lf 2n/3 \rf}{\lf 2k/3 \rf} & \text{else,}\\
\end{cases}
$$
where $\lc\cdot\rc$ and $\lf\cdot\rf$ are the ceiling and floor functions, respectively.
\ele
\bprf
From Theorem~\ref{ss} and the discussion before this lemma, it suffices to show that the right-hand side of the above congruence is the number of fixed points $f(n,k)$ of the action of $\io$ on $\cT_\la\times\cD_{\la^*}$  for all $\la\sbe k\times(n-k)$.  To this end, consider the first quadrant of the plane together with the lines $x=c$ for $c\Cong_3 0$ or $1$ and the lines $y=d$ for $d\Cong_3 0$ or $2$.  Comparing this description with~\ree{fp}, we see that the vertical lines give exactly the row lengths for rows fixed by $\io$ and the horizontal lines do the same for the columns of $\la^*$.  It follows that $f(n,k)$ is just the number of lattice paths from $(0,0)$ to $(n-k,k)$ taking unit steps east or north and staying on the given lines.  Thus $f(n,k)=0$ if $n\Cong_3 0$ and $k\Cong_3 1$ since then $(n-k,k)$ is on none of the lines.  The other cases are handled similarly.
\eprf

\noindent\emph{Proof 1 (of Theorem~\ref{main}).}  This is just a matter of combining Theorem~\ref{binom} and Lemma~\ref{lattice}.  There are three cases, all of them similar, so we will only do the one corresponding to the third expression in the lemma.  If $0\le n, k < 3\cdot 2^m$ for some $m\ge0$ then $0\le   \lf  2n/3 \rf, \lf 2k/3\rf< 2^{m+1}$. Thus
$$
\binom{n+3\cdot  2^m}{k}_F \Cong_2  \binom{\lf 2n/3 \rf + 2^{m+1}}{\lf 2k/3 \rf}\Cong_2 
 \binom{\lf 2n/3 \rf}{\lf 2k/3 \rf} \Cong_2 \binom{n}{k}_F
$$
as desired.\hqed

%
%

\section{An analogue of a congruence of Lucas}
\label{acl}

Our next proof of Theorem~\ref{main} will be number theoretic, based on Fibonomial analogues of two well-known congruences for the ordinary binomial coefficients due to Lucas and Kummer.  To state these, we will need to discuss expansions of integers in various bases.

Suppose $\bb=(b_0,b_1,b_2,\ldots)$ is an infinite increasing sequence of positive integers such that $b_0 = 1$ and for $i\ge1$ we have $b_{i-1}|b_i$.  Then every positive integer $n$ has an \emph{expansion in base $\bb$} which is the unique way of writing
$n=n_0 b_0 + n_1 b_1 + \cdots$ where $0\le n_i<b_{i+1}/b_i$ for all $i$.  
In this case, we write $(n)_\bb=(n_0,n_1,n_2,\ldots)=(n_i)_{i\ge0}$.  If $b_i=m^i$ for some integer $m$ then we use $(n)_m$ for $(n)_\bb$.  In this case, one can extend the expansion to all rational numbers in the usual way.  There is another base which will be useful to us as the reader might expect, namely
$$
\bF=(1,3,3\cdot 2, 3\cdot 2^2, \ldots).
$$

There are two famous theorems due to Kummer and Lucas about congruence properties of binomial coefficients.  To state the former, we need the \emph{$p$-adic valuation} of an integer $n$ which is
$$
\nu_p(n)=\text{the highest power of $p$ dividing $n$.}
$$
\bth[Kummer~\cite{kum:ear}]
If $p$ is  prime then $\nu_p(\binom{m+n}{m})$ is the number of carries in doing the addition $(m)_p+(n)_p$.\hqedm

\eth
\bth[Lucas~\cite{luc:cne}]
Let $p$ be prime and $(n)_p=(n_i)_{i\ge0}$, $(k)_p=(k_i)_{i\ge0}$.  Then

\medskip

\eqqed{
\binom{n}{k} \Cong_p \binom{n_0}{k_0}\binom{n_1}{k_1}\binom{n_2}{k_2}\cdots
}
\eth

Knuth and Wilf generalized Kummer's Theorem to a class of binomial coefficient analogues.  We will need the following particular case of one of their theorems.
\bth[Knuth and Wilf~\cite{kw:ppd}]
\label{kw}
The valuation  $\nu_2(\binom{m+n}{m}_F)$ is the number of carries  in doing the addition $(m/3)_2+(n/3)_2$ where carries to the right of the radix point are not counted and an extra $1$ is added if there is a carry from the one's place to the two's place.\hqed
\eth

We can now prove our analogue of Lucas' Theorem for the prime $p=2$ and the  base $\bF$.
\bth
\label{LucF}
Let $(n)_\bF=(n_i)_{i\ge0}$ and $(k)_\bF=(k_i)_{i\ge0}$.  Then

\medskip

\eqqed{
\binom{n}{k}_F \Cong_2 \binom{n_0}{k_0}_F\binom{n_1}{k_1}_F\binom{n_2}{k_2}_F\cdots
}
\eth
\bprf
It suffices to show that the right-hand side is $0$ or $1$ depending on whether there is a carry or not as described in Theorem~\ref{kw}.  The proof now becomes a case-by-case analysis depending on the congruence classes of $n$ and $k$ modulo $3$.  Since they are all similar, we will just do the case $n\Cong_3 1$ and $k\Cong_3 2$ to illustrate.  But then $n-k\Cong_3 2$ and so there will be a carry across the radix point when doing the addition $(k/3)_2+((n-k)/3)_2$  This implies 
$\binom{n}{k}_F \Cong_2 0$ on the left-hand side.  Considering the product on the right-hand side, we see that it will have a factor of $\binom{n_0}{k_0}_F = \binom{1}{2}_F = 0$.  Thus the two sides agree.
\eprf

\noindent\emph{Proof 2 (of Theorem~\ref{main}).}  Since $0\le n,k < 3\cdot 2^m$ we have $(n)_\bF=(n_0,\ldots,n_m)$ and 
$(k)_\bF=(k_0,\ldots,k_m)$ where final digits equal to zero have been ignored.  It follows that 
$(n+3\cdot 2^m)_\bF= (n_0,\ldots,n_m,1)$.  So applying Theorem~\ref{LucF}
$$
 \binom{n+ 3\cdot 2^m}{k}_F  \Cong_2 
 \binom{n_0}{k_0}_F \cdots \binom{n_m}{k_m}_F \binom{1}{0}_F=
 \binom{n_0}{k_0}_F \cdots \binom{n_m}{k_m}_F \Cong_2
 \binom{n}{k}_F 
$$
which is what we want.\hqed

%
%

\section{An inductive proof}
\label{ip}

In order to give our inductive proof, we need the recursion satisfied by the Fibonomials.  One can find a proof of the following result, e.g., in~\cite{ss:cib}.
\bpr
\label{rr}
We have $\binom{0}{0}_F=1$ and, for $n\ge1$,

\medskip

\eqqed{
\binom{n}{k}_F=F_{n-k+1} \binom{n-1}{k-1}_F + F_{k-1} \binom{n-1}{k}_F.
}
\epr

\noindent\emph{Proof 3 (of Theorem~\ref{main}).}  The base case is trivial, so assume the result for values less than 
$n+3\cdot 2^m$ where $0\le n< 3\cdot 2^m$.  Also suppose $0\le k < 3\cdot 2^m$.  The proof for $n=0$ is only slightly different from the proof for $0<n< 3\cdot 2^m$, so we will just do the latter.  Note that by Proposition~\ref{Fn2}, if the subscripts of two Fibonacci numbers differ by a multiple of three, then the Fibonacci numbers themselves have the same parity.  Using this observation, Proposition~\ref{rr}, and induction gives
\begin{align*}
\binom{n+3\cdot 2^m}{k}_F&=F_{n-k+3\cdot 2^m+1}\binom{n+3\cdot 2^m-1}{k-1}_F + F_{k-1} \binom{n+3\cdot 2^m-1}{k}_F\\[10pt]
&\Cong_2 F_{n-k+1}\binom{n-1}{k-1}_F + F_{k-1} \binom{n-1}{k}_F\\[10pt]
&= \binom{n}{k}_F
\end{align*}
which is the desired conclusion.\hqed

%
%

\bfi
$$
\barr{ccccccccccccccccccccccc}
&&&&&&&&&&&1\\
&&&&&&&&&&1&&1\\
&&&&&&&&&1&&1&&1\\
&&&&&&&&1&&-1&&-1&&1\\
&&&&&&&1&&0&&0&&0&&1\\
&&&&&&1&&-1&&0&&0&&-1&&1\\
&&&&&1&&-1&&1&&0&&1&&-1&&1&\\
&&&&1&&1&&-1&&-1&&-1&&-1&&1&&1\\
&&&1&&0&&0&&0&&-1&&0&&0&&0&&1\\
&&1&&1&&0&&0&&1&&1&&0&&0&&1&&1\\
&1&&1&&1&&0&&-1&&-1&&-1&&0&&1&&1&&1\\
1&&-1&&-1&&1&&1&&-1&&-1&&1&&1&&-1&&-1&&1\\
\rule{8pt}{0pt}&\rule{8pt}{0pt}&\rule{8pt}{0pt}&\rule{8pt}{0pt}&\rule{8pt}{0pt}&\rule{8pt}{0pt}&\rule{8pt}{0pt}&
\rule{8pt}{0pt}&\rule{8pt}{0pt}&\rule{8pt}{0pt}&\rule{8pt}{0pt}&\rule{8pt}{0pt}&\rule{8pt}{0pt}&\rule{8pt}{0pt}&
\rule{8pt}{0pt}&\rule{8pt}{0pt}&\rule{8pt}{0pt}&\rule{8pt}{0pt}&\rule{8pt}{0pt}&\rule{8pt}{0pt}&\rule{8pt}{0pt}&
\rule{8pt}{0pt}&\rule{8pt}{0pt}
\earr
$$
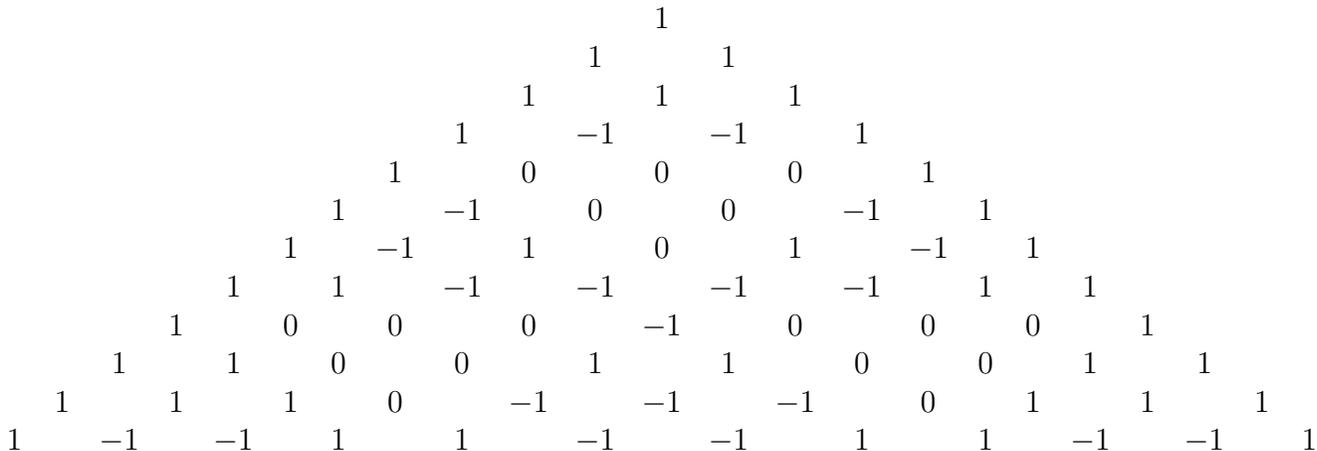
\capt{The Fibonomial triangle  modulo $3$}
\label{tri3} 
\efi

\section{The modulo three case}
\label{mtc}

As mentioned in the introduction, one can extend Theorem~\ref{binom} to Pascal's triangle modulo any prime $p$~\cite{sve:dv}.  If one considers the triangle of the first $p^m$ rows, then the triangle of the first $p^{m+1}$ rows breaks into triangles which are scalar multiples of the smaller triangle with inverted triangles of zeros in between.  

For the the Fibonomial triangle, the situation is more complicated. Here we will consider the case $p=3$ and postpone discussion of general $p$ until the next section.  For $p=3$, the sides of the triangles are of size $4\cdot 3^m$, $m\ge0$.  Figure~\ref{tri3} shows the first $12$ rows in this case.  Note that the triangle $T$ of side $4$ is repeated at the bottom left and bottom right of the triangle of side $12$.  But the other three triangles alternate lines (either rows or diagonals) of $T$ with lines which are negatives of the corresponding line of $T$.  One can prove that this behavior continues by a messy case-by-base induction.  Since the demonstration of each case does not differ significantly from the third proof of Theorem~\ref{main}, we will suppress the details and just state the result.

\bth
\label{mod3}
Given $m\ge0$ and $0\le n,k<4\cdot 3^m$ we have
$$
\binom{n+4\cdot 3^m}{k}_F\Cong_3
\case{\hs{10pt}\dil\binom{n}{k}_F}{if $k$ is even,\rule{0pt}{20pt}}{\dil -\binom{n}{k}_F}{if $k$ is odd,\rule{0pt}{25pt}}
$$
and
$$
\binom{n+8\cdot 3^m}{k}_F\Cong_3 \binom{n}{k}_F,
$$
and
$$
\binom{n+8\cdot 3^m}{k+4\cdot 3^m}_F\Cong_3
\case{\dil -\binom{n}{k}_F}{if $n$ is even,\rule{0pt}{20pt}}{\hs{10pt}\dil \binom{n}{k}_F}{if $n$ is odd.\rule{0pt}{25pt}}
$$
All other values in rows $4\cdot 3^m,\dots,4\cdot 3^{m+1}-1$ are determined by symmetry.\hqed
\eth

%
%

\section{Open questions}
\label{oq}

\subsection{Modulo $3$ redux}

It would be interesting to find proofs of Theorem~\ref{mod3} using combinatorial or number theoretic means.  The combinatorial demonstration  of Theorem~\ref{main} goes through until the one gets to the analogue of Lemma~\ref{lattice}.  The fixed points are now counted by lattice paths where the steps are weighted either $1$ or $2$ and there does not seem to be an easy description of the sum of the weights of paths to a given vertex.

As far as the number theoretic proof, the obvious guess for a modulo $3$ version of Theorem~\ref{LucF} is false.  More specifically, consider the base  $\bT=(1,4,4\cdot 3, 4\cdot 3^2,\dots)$.  If $(n)_\bT=(n_i)_{i\ge0}$ and $(k)_\bT=(k_i)_{i\ge0}$ then one only appears to have
$$
\binom{n}{k}_F \Cong_3 \pm \binom{n_0}{k_0}_F\binom{n_1}{k_1}_F\binom{n_2}{k_2}_F\cdots,
$$
and determining which sign to use does not seem to be an easy matter.  But maybe some other base is called for.

\subsection{Higher modulus}

Having an analogue of Theorem~\ref{main} for any prime modulus $p$ would be quite interesting.  However, this is probably a very difficult problem for the following reason.  It is not hard to show that, given $p$, there must be a least positive integer $p^*$ such that $p$ divides $F_{p^*}$.  We have seen that $2^*=3$ and $3^*=4$.  In fact, $p$ divides $F_n$ if and only if $p^*$ divides $n$; see the paper of Robinson~\cite{rob:fmm} for details.  However, determining $p^*$ for an arbitrary prime is a well-studied open problem.

\subsection{Valuations}

Together with Amdeberhan and Moll~\cite{acms:pgf}, we are considering $2$-adic valuations of generalized Fibonomial coefficients and related sequences.  One can generalize the Fibonacci numbers by considering a sequence of polynomials $\{n\}$, $n\ge0$, in the variables $s$ and $t$ defined by $\{0\}=0$, $\{1\}=1$, and
$$
\{n\}=s\{n-1\}+t\{n-2\}
$$
for $n\ge 2$.  The analogues of the Fibotorials and Fibonomial coefficients are
$$
\{n\}!=\{1\}\{2\}\dots\{n\}
$$
and
$$
\bin{n}{k}=\frac{\{n\}!}{\{k\}!\{n-k\}!},
$$
respectively.  One can show that these are polynomials in $s,t$, and it was in this context that the combinatorial interpretation in~\cite{ss:cib} was given which we have used in  the case $s=t=1$.  We also note that if $s=2$,  $t=-1$ then $\{n\}=n$ and we recover the nonnegative integers.

Louis Shapiro [personal communication] suggested considering the following quotients
$$
C_{\{n\}}=\frac{1}{\{n+1\}}\bin{2n}{n}
$$
which specialize to the Catalan numbers, $C_n$, when $s=2$ and $t=-1$.  It is not hard to show, as done by Ekhad~\cite{ekh:ssl}, that these are always polynomials in $s,t$, thus answering one of the questions raised by Shapiro.

The $C_n$ themselves are known to be odd if and only if $n=2^m-1$ for some $m\ge0$.  In fact, one can express the full $2$-adic valuation using the function
$$
\ze_\bb(n)=\text{the number of nonzero digits in the base $\bb$ expansion of $n$}.
$$
\bth
\label{nuCn}
We have

\eqqed{
\nu_2(C_n)=\ze_2(n+1)-1.
}
\eth

For a (mostly) combinatorial proof of this result, see the paper of Deutsch and Sagan~\cite{ds:ccm}.  
There is an analogue of Theorem~\ref{nuCn}  for $C_{\{n\}}$  using the base $\bF$, as well as related theorems for other values of $s$ and $t$.
\bth[Amdeberhan, Chen, Moll, and Sagan~\cite{acms:pgf}]
Let $s$ be odd and $t\equiv 1 \pmod{8}$.  Then

\medskip

\eqqed{
\nu_2(C_{\{n\}})=\case{\ze_\bF(n+1)}{if $n \Cong 3$ or $4\ (\Mod 6)$,}{\ze_\bF(n+1)-1}{else.}
}
\eth

Shapiro also asked for a combinatorial interpretation of the $C_{\{n\}}$.  Despite the plethora of combinatorial interpretations for the ordinary Catalan numbers, it is still an open problem to give one for general $s$ and $t$.


\begin{thebibliography}{1}

\bibitem{acms:pgf}
Tewodros Amdeberhan, Xi~Chen, Victor Moll, and Bruce~E. Sagan.
\newblock Properties of generalized {F}ibonomial coefficients and related
  sequences.
\newblock In preparation.

\bibitem{ds:ccm}
Emeric Deutsch and Bruce~E. Sagan.
\newblock Congruences for {C}atalan and {M}otzkin numbers and related
  sequences.
\newblock {\em J. Number Theory}, 117(1):191--215, 2006.

\bibitem{ekh:ssl}
Shalosh Ekhad.
\newblock
  The~{S}agan-{S}avage~{L}ucas-{C}atalan~polynomials~have~positive~coefficient%
s.
\newblock Preprint
  {\texttt{http://www.math.rutgers.edu/~zeilberg/mamarim/mamarimhtml/bruce.htm%
l}}.

\bibitem{kw:ppd}
Donald~E. Knuth and Herbert~S. Wilf.
\newblock The power of a prime that divides a generalized binomial coefficient.
\newblock {\em J. Reine Angew. Math.}, 396:212--219, 1989.

\bibitem{kum:ear}
E.~E. Kummer.
\newblock {\"U}ber die {E}rg\"anzungss\"atze zu den allgemeinen
  reciprocit\"atsgesetzen.
\newblock {\em J. reine angew. Math.}, 44:93--146, 1852.

\bibitem{luc:cne}
E.~Lucas.
\newblock Sur les congruences des nombres eul\'eriens et les coefficients
  diff\'erentiels des functions trigonom\'etriques suivant un module premier.
\newblock {\em Bull. Soc. Math. France}, 6:49--54, 1878.

\bibitem{rob:fmm}
D.~W. Robinson.
\newblock The {F}ibonacci matrix modulo {$m$}.
\newblock {\em Fibonacci Quart}, 1(2):29--36, 1963.

\bibitem{ss:cib}
Bruce~E. Sagan and Carla~D. Savage.
\newblock Combinatorial interpretations of binomial coefficient analogues
  related to {L}ucas sequences.
\newblock {\em Integers}, 10:A52, 697--703, 2010.

\bibitem{sve:dv}
Marta Sved.
\newblock Divisibility---with visibility.
\newblock {\em Math. Intelligencer}, 10(2):56--64, 1988.

\end{thebibliography}
\end{document}